\documentclass{article}
\usepackage[english,russian,ukrainian]{babel}
\usepackage{euscript}
\usepackage{latexsym,amssymb,amsmath}
\usepackage[cp1251]{inputenc}
\usepackage{amsfonts,amssymb}
\usepackage{graphicx,graphics,hhline}
\usepackage{euscript}
\newcounter{theorem} 
\newcounter{lemma} 

\,
\,

\begin{document}

\vspace*{7mm}

\LARGE

\addtocounter{page}{-1}
\thispagestyle{empty}

\begin{center}
\textbf{Integral theorems for monogenic functions in
commutative algebras}
\end{center}
\vskip4mm
\large
\begin{center}\textbf{V.\,S.\,Shpakivskyi}\end{center}

\vspace{7mm}
\small
Let $\mathbb{A}_n^m$ be an arbitrary $n$-dimensional
commutative associative algebra
over the field of complex numbers with $m$ idempotents. Let
$e_1=1,e_2,\ldots,e_k$ with $2\leq k\leq 2n$ be elements of $\mathbb{A}_n^m$ which are linearly
independent over the field of real numbers. We consider monogenic
(i.~e. continuous and differentiable in the sense of Gateaux)
functions of the variable $\sum_{j=1}^k x_j\,e_j$\,, where $x_1,x_2,\ldots,x_k$ are
real, and we prove curvilinear analogues of the
Cauchy integral theorem, the Morera theorem and the Cauchy integral formula in $k$-dimensional ($2\leq k\leq 2n$)  real subset of the algebra $\mathbb{A}_n^m$. The present article is generalized of the author's paper [\ref{Shpakivskyi-AACA-2015}], where mentioned results are obtained for $k=3$.
\vspace{7mm}

\large

\section{Introduction}
\vskip2mm

The Cauchy integral theorem and Cauchy integral
formula for the holomorphic function of the complex variable are
a fundamental result of the classical complex analysis.
Analogues of these results are also an important tool in commutative
algebras of dimensional more that $2$.

In the paper of E.~R.~Lorch [\ref{Lorch}] for functions differentiable in the
sense of Lorch in an arbitrary convex domain of commutative
associative Banach algebra, some properties similar to properties
of holomorphic functions of complex variable (in particular, the
curvilinear integral Cauchy theorem and the integral Cauchy formula, the
Taylor expansion and the Morera theorem) are established. E.~K.~Blum [\ref{Blum}]
withdrew a convexity condition of a domain in the
mentioned results from [\ref{Lorch}].

Let us note that {\it a
priori}\/ the differentiability of a function in the
sense of Gateaux is a restriction weaker than the
differentiability of this function in the sense of Lorch.

Therefore, we consider a \textit{monogenic} functions defined as a continuous
and differentiable in the sense of Gateaux. Also we assume that a
monogenic function is given in a domain of three-dimensional subspace of
an arbitrary commutative associative algebra with unit over the field of complex
 numbers.
 In this situation
the results established in the papers [\ref{Lorch}, \ref{Blum}] is not applicable
 for a mentioned monogenic function, because it deals with an
 integration along a curve on which the function is not given, generally speaking.

 In the papers [\ref{Pl-Shp3}, \ref{Pl-Shp-Algeria}, \ref{Pl-Pukh-Analele}] for monogenic
function the curvilinear analogues of the
Cauchy integral theorem, the Cauchy integral formula and the
 Morera theorem are obtained in special finite-dimensional
 commutative associative algebras.  The results of the papers [\ref{Pl-Shp3}, \ref{Pl-Shp-Algeria}, \ref{Pl-Pukh-Analele}] are generalized in the article [\ref{Shpakivskyi-AACA-2015}] to an arbitrary commutative associative algebra. At the same time in [\ref{Shpakivskyi-AACA-2015}] monogenic functions were defined in a domain of real three-dimensional subspace  of an algebra.

 In this paper we generalize results of the papers
[\ref{Shpakivskyi-AACA-2015}] assuming that monogenic functions are defined in a domain of real $k$-dimensional subspace  of an algebra.

 Let us note that some analogues of the curvilinear  Cauchy's integral theorem
 and the Cauchy's integral formula for another classes of functions in
 special commutative algebras are
established in the papers [\ref{Ketchum-28}, \ref{Ketchum-29}, \ref{Rosculet-54}, \ref{Rosculet-55}].\vskip2mm

\section{The algebra $\mathbb{A}_n^m$}
\vskip2mm

Let $\mathbb{N}$ be the set of natural numbers.
We fix the numbers $m,n\in\mathbb{N}$ such that $m\leq n$.
Let $\mathbb{A}_n^m$ be an arbitrary commutative associative algebra with
 unit over the field of complex number $\mathbb{C}$.
   E.~Cartan [\ref{Cartan}, pp.~33 -- 34]
  proved that in the algebra $\mathbb{A}_n^m$ there exist a basis $\{I_k\}_{k=1}^{n}$
  satisfies the following multiplication rules:
  \vskip3mm
1.  \,\,  $\forall$\; $r,s\in[1,m]\cap\mathbb{N}$\,: \qquad $I_rI_s=\left\{
\begin{array}{rcl}
0 &\mbox{if} & r\neq s,\vspace*{2mm} \\
I_r &\mbox{if} & r=s;\\
\end{array}
\right.$

\vskip5mm

2. \,\,  $\forall$\; $r,s\in[m+1,n]\cap\mathbb{N}$\,: \qquad $I_rI_s=
\sum\limits_{k=\max\{r,s\}+1}^n\Upsilon_{r,k}^{s}I_k$\,;

\vskip5mm

3.\,\, $\forall$\; $s\in[m+1,n]\cap\mathbb{N}$\;  $\exists!\;
 u_s\in[1,m]\cap\mathbb{N}$ \;$\forall$\,
 $r\in[1,m]\cap\mathbb{N}$\,:\;\;

\begin{equation}\label{mult_rule_3}
I_rI_s=\left\{
\begin{array}{ccl}
0 \;\;\mbox{if}\;\;  r\neq u_s\,,\vspace*{2mm}\\
I_s\;\;\mbox{if}\;\;  r= u_s\,. \\
\end{array}
\right.\medskip
\end{equation}
Furthermore, the structure constants $\Upsilon_{r,k}^{s}\in\mathbb{C}$
 satisfy the associativity conditions:
\vskip2mm
(A\,1).\,\, $(I_rI_s)I_p=I_r(I_sI_p)$ \; $\forall$\, $r,s,p\in[m+1,n]\cap\mathbb{N}$;
\vskip2mm
(A\,2).\,\, $(I_uI_s)I_p=I_u(I_sI_p)$ \; $\forall$\, $u\in[1,m]\cap\mathbb{N}$\;
 $\forall$\,
 $s,p\in[m+1,n]\cap\mathbb{N}$.
\vskip2mm
Obviously, the first $m$ basis vectors $\{I_u\}_{u=1}^m$ are the
idempotents and, respectively,
form the semi-simple subalgebra. Also the vectors $\{I_r\}_{r=m+1}^n$
form the nilpotent subalgebra of algebra
 $\mathbb{A}_n^m$.
The unit of $\mathbb{A}_n^m$ is the element $1=\sum_{u=1}^mI_u$. Therefore,
we will write that the algebra $\mathbb{A}_n^m$ is a semi-direct sum of the
$m$-dimensional semi-simple subalgebra $S$ and $(n-m)$-dimensional
nilpotent subalgebra $N$, i.~e.
$$\mathbb{A}_n^m=S\oplus_s N.
$$

Let us note that nilpotent algebras are fully described for the dimensions
$1,2,3$ in the paper [\ref{Burde_de_Graaf}], and some four-dimensional nilpotent algebras
can be found in the papers [\ref{Burde_Fialowski}],  [\ref{Martin}].

The algebra $\mathbb{A}_n^m$ contains $m$ maximal ideals
$$\mathcal{I}_u:=\Biggr\{\sum\limits_{k=1,\,k\neq u}^n\lambda_kI_k:\lambda_k\in
\mathbb{C}\Biggr\}, \quad  u=1,2,\ldots,m,
$$
the intersection of which is the radical $$\mathcal{R}:=
\Bigr\{\sum\limits_{k=m+1}^n\lambda_kI_k:\lambda_k\in
\mathbb{C}\Bigr\}.$$

We define $m$ linear functionals $f_u:\mathbb{A}_n^m\rightarrow\mathbb{C}$ by put
$$f_u(I_u)=1,\quad f_u(\omega)=0\quad\forall\,\omega\in\mathcal{I}_u\,,
\quad u=1,2,\ldots,m.
$$
Since the kernels of functionals $f_u$ are, respectively, the maximal ideals
 $\mathcal{I}_u$, then these functionals are also continuous and multiplicative
  (see [\ref{Hil_Filips}, p. 147]).

\vskip2mm
\section{Monogenic functions in $E_k$}
\vskip2mm

Let us consider the vectors $e_1=1,e_2,\ldots,e_k$ in $\mathbb{A}_n^m$, where $2\leq k\leq 2n$, and these vectors are linearly independent over the field of real numbers
$\mathbb{R}$ (see [\ref{Pl-Pukh-Analele}]). It means that the equality
$$\sum\limits_{j=1}^k\alpha_je_j=0,\qquad \alpha_j\in\mathbb{R},$$
holds if and only if $\alpha_j=0$ for all $j=1,2,\ldots,k$.

Let the vectors $e_1=1,e_2,\ldots,e_k$ have the following decompositions with respect to
the basis $\{I_r\}_{r=1}^n$:
\begin{equation}\label{e_1_e_2_e_3-k}
e_1=\sum\limits_{r=1}^mI_r\,, 
\quad e_j=\sum\limits_{r=1}^na_{jr}\,I_r\,,\quad a_{jr}\in\mathbb{C},\quad j=2,3,\ldots,k.
\end{equation}

Let $\zeta:=\sum\limits_{j=1}^kx_j\,e_j$, where $x_j\in\mathbb{R}$. It is
obvious that
 $$\xi_u:=f_u(\zeta)=x_1+\sum\limits_{j=2}^kx_j\,a_{ju},\quad u=1,2,\ldots,m.$$
 Let
 $E_k:=\{\zeta=\sum\limits_{j=1}^kx_je_j:\,\, x_j\in\mathbb{R}\}$ be the
linear span of vectors $e_1=1,e_2,\ldots,e_k$ over the field
$\mathbb{R}$. We note that in the further investigations,
 it is essential assumption:
$f_u(E_k)=\mathbb{C}$ for all\, $u=1,2,\ldots,m$.
 Obviously, it holds if and only if for every fixed $u=1,2, \ldots, m$
at least one of the numbers $a_{2u}$, $a_{3u},\ldots,a_{ku}$ belongs to
$\mathbb{C}\setminus\mathbb{R}$.

With a set $Q_{\mathbb{R}}\subset\mathbb{R}^k$ we associate the set $Q:=
\{\zeta=\sum\limits_{j=1}^kx_je_j:(x_1,\ldots,x_k)\in Q_{\mathbb{R}}\}$ in $E_k$. We also note that the topological properties of a set $Q$ in $E_k$ understood as a
corresponding topological properties of a set $Q_{\mathbb{R}}$ in $\mathbb{R}^k$.
For example, a homotopicity of a curve $\gamma\subset E_k$ to the zero means a
homotopicity of $\gamma_{\mathbb{R}}\subset\mathbb{R}^k$ to the zero; a rectifiability of a curve $\gamma\subset E_k$ we understand as a rectifiability of the curve $\gamma_{\mathbb{R}}\subset\mathbb{R}^k$,  etc.

Let $\Omega$ be a domain in $E_k$. With a domain $\Omega\subset E_k$ we associate the domain  $$\Omega_{\mathbb{R}}:=\Big\{(x_1,x_2,\ldots,x_k)\in\mathbb{R}^k:\,\zeta=\sum\limits_{j=1}^kx_j\,e_j\in\Omega\Big\}$$ in $\mathbb{R}^k$.

We say that a continuous function
$\Phi:\Omega\rightarrow\mathbb{A}_n^m$ is \textit{monogenic}
in $\Omega$ if $\Phi$ is differentiable in the sense of
Gateaux in every point of $\Omega$, i.~e. if  for every
$\zeta\in\Omega$ there exists an element
$\Phi'(\zeta)\in\mathbb{A}_n^m$ such that
\begin{equation}\label{monogennaOZNA}\medskip
\lim\limits_{\varepsilon\rightarrow 0+0}
\left(\Phi(\zeta+\varepsilon
h)-\Phi(\zeta)\right)\varepsilon^{-1}= h\Phi'(\zeta)\quad\forall\,
h\in E_k.\medskip
\end{equation}
$\Phi'(\zeta)$ is the \textit{Gateaux derivative} of the function
$\Phi$ in the point $\zeta$.

Consider the decomposition of a function
$\Phi:\Omega\rightarrow\mathbb{A}_n^m$ with respect to the
basis $\{I_r\}_{r=1}^n$:
\begin{equation}\label{rozklad-Phi-v-bazysi-k}
\Phi(\zeta)=\sum_{r=1}^n U_r(x_1,x_2,\ldots,x_k)\,I_r\,.
 \end{equation}

In the case where the functions $U_r:\Omega_{\mathbb{R}}\rightarrow\mathbb{C}$ are
$\mathbb{R}$-differentiable in $\Omega_{\mathbb{R}}$, i.~e. for every $(x_1,x_2,\ldots,x_k)\in\Omega_{\mathbb{R}}$,
$$U_r\left(x_1+\Delta x_1,x_2+\Delta x_2,\ldots,x_k+\Delta x_k\right)-U_r(x_1,x_2,\ldots,x_k)=
$$
$$=\sum\limits_{j=1}^k\frac{\partial U_r}{\partial x_j}\,\Delta x_j+
\,o\left(\sqrt{\sum\limits_{j=1}^k(\Delta x_j)^2}\,\right), \qquad \sum\limits_{j=1}^k(\Delta x_j)^2\to 0\,,$$
the function $\Phi$ is monogenic in the domain $\Omega$ if
and only if the following Cauchy~-- Riemann conditions are
satisfied in $\Omega$:
\begin{equation}\label{Umovy_K-R-k}
\frac{\partial \Phi}{\partial x_j}=\frac{\partial \Phi}{\partial
x_1}\,e_j\qquad\text{for all}\quad j=2,3,\ldots,k.
\end{equation}

An expansion of the resolvent is of the form (see [\ref{Shpakivskyi-Zb-2015-1}]):
 \begin{equation}\label{lem-rez-}
(te_1-\zeta)^{-1}=\sum\limits_{u=1}^m\frac{1}{t-\xi_u}\,I_u+
 \sum\limits_{s=m+1}^{n}\sum\limits_{r=2}^{s-m+1}\frac{Q_{r,s}}
 {\left(t-\xi_{u_{s}}\right)^r}\,I_{s}\,
 \end{equation}
$$\forall\,t\in\mathbb{C}:\,
t\neq \xi_u,\quad u=1,2,\ldots,m,$$
where the coefficients\ $Q_{r,s}$ are determined by the following recurrence
relations:
\begin{equation}\label{Q-}
\begin{array}{c}
\displaystyle
Q_{2,s}=T_s\,,\quad
Q_{r,s}=\sum\limits_{q=r+m-2}^{s-1}Q_{r-1,q}\,B_{q,\,s}\,,\; \;\;r=3,4,\ldots,s-m+1,\\
\end{array}
\end{equation}
with
\begin{equation}\label{B-}
T_s:=\sum\limits_{j=2}^kx_ja_{js}\,, \quad
B_{q,s}:=\sum\limits_{p=m+1}^{s-1}T_p \Upsilon_{q,s}^p\,,
\;\;p=m+2,m+3,\ldots,n,
\end{equation}
 and the natural numbers $u_s$ are defined in the
rule 3 of the multiplication table of algebra
$\mathbb{A}_n^m$.

In the paper [\ref{Shpakivskyi-2014}] an expansion of the resolvent is obtained for $k=3$.

It follows from the relation (\ref{lem-rez-}) that the points
 $(x_1,x_2,\ldots,x_k)\in\mathbb{R}^k$
corresponding to the noninvertible elements
$\zeta=\sum\limits_{j=1}^kx_j\,e_j$
form the set
  \[M_u^{\mathbb{R}}:\quad\left\{
\begin{array}{r}x_1+\sum\limits_{j=2}^kx_j\,{\rm Re}\,a_{ju}=0,\vspace*{3mm} \\
\sum\limits_{j=2}^kx_j\,{\rm Im}\,a_{ju}=0, \\ 
\end{array} \right. \qquad u=1,2,\ldots,m\]
in the $k$-dimensional space $\mathbb{R}^k$.
Also we consider the set
$M_u:=\{\zeta\in E_k: f_u(\zeta)=0\}$ for $u=1,2,\ldots,m$.
It is obvious that the set $M_u^{\mathbb{R}}\subset\mathbb{R}^k$ is congruent with the set  $M_u\subset E_k$.

Denote by $D_u\subset\mathbb{C}$ the image of $\Omega$ under the mapping
$f_u$,\, $u=1,2, \ldots, m$. We say that a domain $\Omega\subset E_k$ is \textit{convex with respect to the set of directions}
$M_u$ if  $\Omega$ contains the segment $\{\zeta_1+\alpha(\zeta_2-\zeta_1):\alpha\in[0,1]\}$ for all $\zeta_1,\zeta_2\in \Omega$ such that $\zeta_2-\zeta_1\in M_u$.
A constructive description of all monogenic functions in the algebra
$\mathbb{A}_n^m$ by means of holomorphic functions of the complex variable are
  obtained in the paper [\ref{Shpakivskyi-Zb-2015-1}].  Namely, it is
 proved the theorem:

Let a domain $\Omega\subset E_k$ be convex with respect to the set of directions
$M_u$  and $f_u(E_k)=\mathbb{C}$ for all
 $u=1,2,\ldots, m$.  Then every
monogenic function $\Phi:\Omega\rightarrow\mathbb{A}_n^m$
can be expressed in the form
 \begin{equation}\label{Teor--1-k-}
\Phi(\zeta)=\sum\limits_{u=1}^mI_u\,\frac{1}{2\pi
i}\int\limits_{\Gamma_u} F_u(t)(te_1-\zeta)^{-1}\,dt+
\sum\limits_{s=m+1}^nI_s\,\frac{1}{2\pi i}\int\limits_
{\Gamma_{u_s}}G_s(t)(te_1-\zeta)^{-1}\,dt,
 \end{equation}\vskip1mm
\noindent where $F_u$ and $G_s$ are certain holomorphic functions in the
domains $D_u$ and $D_{u_s}$, respectively, and $\Gamma_q$ is a
closed Jordan rectifiable curve in $D_q$ which surrounds the point
$\xi_q$ and contains no points $\xi_{\ell}$, $\ell,q=1,2,\ldots,
m$,\,$\ell\neq q$.

We note that in the paper [\ref{Shpakivskyi-2014}] the previous result is proved for $k=3$.

\vskip2mm
\section{Cauchy integral theorem for a curvilinear integral}
\vskip2mm

 Let $\gamma$ be a Jordan rectifiable curve in $E_k$. For a
continuous function $\Psi:\gamma\rightarrow
\mathbb{A}_n^m$ of the form
\begin{equation}\label{Phi-form-}
\Psi(\zeta)=\sum\limits_{r=1}^{n}U_r(x_1,x_2,\ldots,x_k)\,I_r+i\sum
\limits_{r=1}^{n}V_r(x_1,x_2,\ldots,x_k)\,I_r\,,
\end{equation}
where $(x_1,x_2,\ldots,x_k)\in\gamma_{\mathbb{R}}$ and $U_r:\gamma_{\mathbb{R}}\rightarrow\mathbb{R}$,
$V_r:\gamma_{\mathbb{R}}\rightarrow\mathbb{R}$,
we define an integral along a Jordan rectifiable curve $\gamma$ by
 the equality:
$$\int\limits_{\gamma}\Psi(\zeta)d\zeta:=\sum\limits_{j=1}^{k}e_j\sum\limits_{r=1}^{n}
I_{r}\int\limits_{\gamma_{\mathbb{R}}}U_r(x_1,x_2,\ldots,x_k)dx_j+$$
$$+
i\sum\limits_{j=1}^{k}e_j\sum\limits_{r=1}^nI_{r}\int\limits_{\gamma_{\mathbb{R}}}V_r
(x_1,x_2,\ldots,x_k)dx_j\,,$$
where $d\zeta:=e_1dx_1+e_2dx_2+\ldots+e_kdx_k$.

Also we define a surface integral. Let $\Sigma$ be a piece-smooth hypersurface in
 $E_k$. For a continuous
function $\Psi:\Sigma\rightarrow \mathbb{A}_n^m$ of the
form (\ref{Phi-form-}), where $(x_1,x_2,\ldots,x_k)\in\Sigma_{\mathbb{R}}$ and $U_r:\Sigma_{\mathbb{R}}\rightarrow\mathbb{R}$,
$V_r:\Sigma_{\mathbb{R}}\rightarrow\mathbb{R}$, we define a surface
integral on $\Sigma$ with the differential form
$dx_p\wedge dx_q$, by the equality
$$\int\limits_{\Sigma}\Psi(\zeta)dx_p\wedge dx_q:=
\sum\limits_{r=1}^nI_{r}\int\limits_{\Sigma_{\mathbb{R}}}U_{r}(x_1,x_2,\ldots,x_k)dx_p\wedge dx_q+$$
$$+
i\sum\limits_{r=1}^nI_{r}\int\limits_{\Sigma_{\mathbb{R}}}V_{r}(x_1,x_2,\ldots,x_k)dx_p\wedge dx_q.
$$

If a function $\Phi:\Omega\rightarrow\mathbb{A}_n^m$ is
continuous together with partial derivatives of the first order in
a domain $\Omega$, and $\Sigma$ is a piece-smooth hypersurface
in $\Omega$, and the edge $\gamma$ of surface $\Sigma$ is a
rectifiable Jordan curve, then the following analogue of the
Stokes formula is true:
$$\int\limits_{\gamma}\Psi(\zeta)d\zeta=\int\limits_
{\Sigma}\left(\frac{\partial\Psi}{\partial
x_1}\,e_{2}-\frac{\partial\Psi}{\partial x_2}\,e_1\right)dx_1\wedge dx_2
+$$
\begin{equation}\label{form-Stoksa-}+
\left(\frac{\partial\Psi}{\partial
x_2}e_{3}-\frac{\partial\Psi}{\partial x_3}e_{2}\right)dx_2\wedge dx_3+\ldots+\left(\frac{\partial\Psi}{\partial
x_k}e_1-\frac{\partial\Psi}{\partial x_1}e_k\right)dx_k\wedge dx_1.
\end{equation}

Now, the next theorem is a result of the formula
(\ref{form-Stoksa-}) and the equalities (\ref{Umovy_K-R-k}).
\vskip2mm

\textbf{Theorem 1.}\label{teo-int-po-kryv-z-neper-poh-}
\emph{ Suppose that $\Phi:\Omega\rightarrow\mathbb{A}_n^m$ is a
monogenic function in a domain $\Omega$, and $\Sigma$ is a
piece-smooth surface in $\Omega$, and the edge $\gamma$ of surface
$\Sigma$ is a rectifiable Jordan curve. Then}
\begin{equation}\label{form-Koshi-po-kryv-}
\int\limits_{\gamma}\Phi(\zeta)d\zeta=0.
\end{equation}

In the case where a domain $\Omega$ is convex, then by the usual way
 (see, e.~g., [\ref{Privalov}]) the equality (\ref{form-Koshi-po-kryv-})
  can be prove for an arbitrary closed Jordan rectifiable curve $\gamma_\zeta$.

In the case where a domain $\Omega$ is an arbitrary, then similarly
 to the proof of Theorem 3.2 [\ref{Blum}] we can
prove the following\vskip2mm

\textbf{Theorem 2.}\label{teo-int-po-kryv-Blum-}
 \emph{ Let $\Phi:\Omega\rightarrow\mathbb{A}_n^m$ be a monogenic
function in a domain $\Omega$. Then for every closed
Jordan rectifiable curve $\gamma$ homotopic to a point in
$\Omega$, the equality \em (\ref{form-Koshi-po-kryv-}) \em is true.}

\vskip2mm
\section{The Morera theorem}
\vskip2mm

To prove the analogue of Morera theorem in the algebra $\mathbb{A}_n^m$,
we introduce auxiliary notions and prove some auxiliary statements.

Let us consider the algebra $\mathbb{A}_n^m(\mathbb{R})$ with  the
 basis $\{I_k,iI_k\}_{k=1}^n$ over the field $\mathbb{R}$ which is isomorphic to the
algebra $\mathbb{A}_n^m$ over the field $\mathbb{C}$. In the algebra
$\mathbb{A}_n^m(\mathbb{R})$ there exist another basis $\{e_r\}_{r=1}^{2n}$,
where the vectors $e_1,e_2,\ldots,e_k$ are the same as in the Section 3.

For the element $a:=\sum\limits_{r=1}^{2n}a_re_r$,\, $a_r\in\mathbb{R}$ we define
the Euclidian norm $$\|a\|:=\sqrt{\sum\limits_{r=1}^{2n}a_r^2}\,.$$ Accordingly,
$\|\zeta\|=\sqrt{\sum\limits_{j=1}^kx_j^2}$\, and $\|e_j\|=1$ for all $j=1,2,\ldots,k$.

Using the Theorem on equivalents of norms, for the element
$b:=\sum\limits_{r=1}^{n}(b_{1r}+ib_{2r})I_r$,\, $b_{1r},b_{2r}\in\mathbb{R}$ we have
the following inequalities
\begin{equation}\label{ner-dlja-integrala-dop-}
|b_{1r}+ib_{2r}|\leq\sqrt{\sum\limits_{r=1}^{2n}\big(b_{1r}^2+b_{2r}^2\big
)}\,\leq c \|b\|,
\end{equation}
where $c$ is a positive constant does not depend on $b$.
\vskip2mm

\textbf{Lemma 1.}\label{lem-22-} \textit{If $\gamma$ is a closed Jordan
rectifiable curve in $E_k$ and function
$\Psi:\gamma\rightarrow
\mathbb{A}_n^m$ is continuous, then}

\begin{equation}\label{ner-dlja-integrala-}
\Biggr\|\int\limits_{\gamma}\Psi(\zeta)\,d\zeta\Biggr\|\leq
c \int\limits_{\gamma}\|\Psi(\zeta)\|\,\|d\zeta\|,
\end{equation}
\textit{where
$c$ is a positive absolutely constant.}

\vskip1mm
\textbf{\textit{Proof.}} Using the representation of function $\Psi$ in the form
 (\ref{Phi-form-}) for $(x_1,x_2,\ldots,x_k)\in\gamma$, we obtain \vspace{2mm}
 $$\Biggl\|\int\limits_{\gamma}\Psi(\zeta)d\zeta\Biggr\|\le
 \sum\limits_{r=1}^{n}\|e_1I_r\|\int\limits_{\gamma_{\mathbb{R}}}\bigl|U_r(x_1,x_2,\ldots,x_k)+
 iV_r(x_1,x_2,\ldots,x_k)\bigr|
 \,dx_1+$$\vspace{1mm}
$$\ldots+\sum\limits_{r=1}^{n}\|e_kI_r\|\int\limits_{\gamma_{\mathbb{R}}}\bigl|U_r(x_1,x_2,\ldots,x_k)+
 i\,V_r(x_1,x_2,\ldots,x_k)\bigr|\,dx_k.$$
Now, taking into account the inequality (\ref{ner-dlja-integrala-dop-}) for
 $b=\Psi(\zeta)$ and the inequalities $\|e_jI_r\|\leq c_j$,\, $j=1,2,\ldots,k$,
 where
$c_j$ are positive absolutely constants,
 we obtain the relation (\ref{ner-dlja-integrala-}).
\noindent The lemma is proved.

We understand a triangle $\triangle$ as a plane figure bounded by three line segments connecting three its vertices. Denote by $\partial\triangle$ the boundary of triangle $\triangle$ in relative topology of its plane.

Using Lemma 1, for functions taking values in the algebra $\mathbb{A}_n^m$, the
following Morera theorem can be established in the usual way.
\vskip2mm

\textbf{Theorem 3.}\label{teo_Morera}
 \emph{ If a function $\Phi:\Omega\rightarrow\mathbb{A}_n^m$ is
continuous in a domain $\Omega$ and satisfies the equality
\begin{equation} \label{Morera}
\int\limits_{\partial\triangle}\Phi(\zeta)d\zeta=0
\end{equation}
 for every triangle $\triangle$ such that closure
$\overline{\triangle}\subset\Omega$, then the function
$\Phi$ is monogenic in the domain $\Omega$.}

\vskip2mm
\section{Cauchy integral formula for a curvilinear integral}
\vskip2mm

 Let $\zeta_0:=\sum\limits_{j=1}^kx_j^{(0)}\,e_j$ be a point in a domain
$\Omega\subset E_k$. Let us take an any 2-dimensional plan containing the point $\zeta_0$ and in this plane we take a circle $C(\zeta_0,\varepsilon)$ of radius $\varepsilon$ with
the center at the point $\zeta_0$, such that this circle completely contained in
$\Omega$. By $C_u(\xi^{(0)}_u,\varepsilon)\subset\mathbb{C}$ we denote the image of $C(\zeta_0,\varepsilon)$
under the mapping $f_u$, $u=1,2,\ldots,m$.  We assume that
the circle $C(\zeta_0,\varepsilon)$ \emph{embraces the set} $\{\zeta-\zeta_0:
\zeta\in\bigcup\limits_{u=1}^m M_u\}$.
It means that the curve $C_u(\xi^{(0)}_u,\varepsilon)$ bounds some domain $D_u'$ and
$f_u(\zeta_0)=\xi^{(0)}_{u}\in D_u'$, \, $u=1,2,\ldots,m$.

We say that the curve $\gamma\subset\Omega$ \emph{embraces once the set} $\{\zeta-\zeta_0:\zeta\in\bigcup\limits_{u=1}^m M_u\}$, if there exists a
 circle $C(\zeta_0,\varepsilon)$
which embraces the mentioned set and is homotopic to $\gamma$ in the domain $\Omega\setminus\{\zeta-\zeta_0:\zeta\in\bigcup\limits_{u=1}^m M_u\}$.

Since the function $\zeta^{-1}$ is continuous on the curve $C(0,\varepsilon)$, then there exists the integral
\begin{equation}\label{lambda-}
\lambda:=\int\limits_{C(0,\varepsilon)}\zeta^{-1}d\zeta.
\end{equation}

The following theorem is an analogue of Cauchy integral
theorem for monogenic function
 $\Phi:\Omega\rightarrow\mathbb{A}_{n}^m$.\vskip2mm

\textbf{Theorem 4.}\label{teo-formula-Koshi-} \textit{Suppose that a domain $\Omega\subset
E_k$ is convex with respect to the set of directions
$M_u$  and $f_u(E_k)=\mathbb{C}$ for all
 $u=1,2,\ldots, m$. Suppose also that
 $\Phi:\Omega\rightarrow\mathbb{A}_n^m$ is a monogenic
function in $\Omega$.  Then for every point
$\zeta_{0}\in\Omega$ the following equality is true:
\begin{equation}\label{form-Koshi-}
\lambda\,\Phi(\zeta_{0})=
\int\limits_{\gamma}\Phi(\zeta)\left(\zeta-\zeta_{0}\right)^{-1}d\zeta,
\end{equation}
where $\gamma$ is an arbitrary closed Jordan rectifiable curve in
$\Omega$, that embraces once the set $\{\zeta-\zeta_0: \zeta\in
\bigcup\limits_{u=1}^m M_u\}$.}
\vskip 1mm

The proof is similar to the proof of Theorem 4 of [\ref{Shpakivskyi-AACA-2015}].

\vskip2mm
\section{A constant $\lambda$}
\vskip2mm

In some special algebras (see [\ref{Pl-Shp3}, \ref{Pl-Shp-Algeria}, \ref{Pl-Pukh-Analele}])
the Cauchy integral formula (\ref{form-Koshi-}) has the form
\begin{equation}\label{form-Koshi--}
\Phi(\zeta_{0})=\frac{1}{2\pi i}
\int\limits_{\gamma_{\zeta}}\Phi(\zeta)\left(\zeta-\zeta_{0}\right)^{-1}d\zeta,
\end{equation}
i.~e.
\begin{equation}\label{lambda-}
\lambda=2\pi i.
\end{equation}

In this Section we indicate a set of algebras $\mathbb{A}_n^m$ for which
(\ref{lambda-}) holds. In this a way we first consider some auxiliary statements.

As a consequence of the expansion (\ref{lem-rez-}), we obtain the
following equality:
\begin{equation}\label{obr-elem-}
\zeta^{-1}=\sum\limits_{r=1}^n\widetilde{A}_r\,I_r
\end{equation}
with the coefficients\, $\widetilde{A}_r$\, determined by the following
 relations:

\begin{equation}\label{A__p-}
\begin{array}{c}
\displaystyle
\widetilde{A}_u=\frac{1}{\xi_u}\,,\;\;u=1,2,\ldots,m, \quad
 \vspace*{4mm}\\
 \displaystyle
\widetilde{A}_s=\sum\limits_{r=2}^{s-m+1}\frac{\widetilde{Q}_{r,s}}{\xi_{u_s}^r}\,,\quad s=m+1,m+2,\ldots,n,
\end{array}
\end{equation}
where $\widetilde{Q}_{r,s}$ are determined by the following recurrence
relations:
\begin{equation}\label{lem_3_2---}
\begin{array}{c}
\displaystyle
\widetilde{Q}_{2,s}=-T_s\,,\quad
\widetilde{Q}_{r,s}=-\sum\limits_{q=r+m-2}^{s-1}\widetilde{Q}_{r-1,q}\,B_{q,\,s}\,,\; \;\;r=3,4,\ldots,s-m+1,\\
\end{array}
\end{equation}
where $T_s$ and $B_{q,s}$ are the same as in the equalities (\ref{B-}),
 and natural numbers $u_s$ are defined in the
rule  3  of the multiplication table of the algebra $\mathbb{A}_n^m$.

Taking into account the equality
(\ref{obr-elem-}) and the relation
$$d\zeta=\sum\limits_{j=1}^kdx_j\,e_j=\sum\limits_{u=1}^m\Big(dx_1+\sum\limits_{j=2}^kdx_j\,a_{ju} \Big)I_u+$$
$$+\sum\limits_{r=m+1}^n\sum\limits_{j=2}^kdx_j\,a_{js}\,I_r
=\sum\limits_{u=1}^md\xi_u\,I_u+\sum\limits_{r=m+1}^ndT_r\,I_r\,,
$$
we have the following equality
$$
\zeta^{-1}d\zeta=\sum\limits_{u=1}^m\widetilde{A}_u\,d\xi_u\,I_u+
\sum\limits_{r=m+1}^n\widetilde{A}_{u_r}\,dT_r\,I_r+$$
\begin{equation}\label{obr-elem-0}
+\sum\limits_{s=m+1}^n\widetilde{A}_{s}\,d\xi_{u_s}\,I_s+
\sum\limits_{s=m+1}^n\sum\limits_{r=m+1}^n\widetilde{A}_{s}\,dT_r\,I_sI_r=:
\sum\limits_{r=1}^n\sigma_r\,I_r\,.
\end{equation}\vskip 2mm

Now, taking into account the denotation (\ref{obr-elem-0}) and the equality
(\ref{A__p-}), we calculate:
$$\int\limits_{C(0,R)}\sum\limits_{u=1}^m\sigma_u\,I_u=
\sum\limits_{u=1}^mI_u\int\limits_{C_u(\xi_u,R)}\frac{d\xi_u}{\xi_u}=
2\pi i\sum\limits_{u=1}^mI_u=2\pi i.
$$

Thus,
\begin{equation}\label{obr-elem--1}
\lambda=2\pi i+\sum\limits_{r=m+1}^nI_r\int\limits_{C(0,R)}\sigma_r\,.
\end{equation}
Therefore, the equality (\ref{lambda-}) holds if and only if
\begin{equation}\label{obr-elem-1}
\int\limits_{C(0,R)}\sigma_r=0\qquad \forall\;r=m+1,\ldots,n.
\end{equation}
  But,
for satisfying the equality (\ref{obr-elem-1}) the differential form $\sigma_r$ must be a total differential of some function.
We note that the property of being a total differential
is invariant under admissible transformations of coordinates
[\ref{Shabat}, p. 328, Theorem 2]. In our situation, if we show that
$\sigma_r$ is a total differential of some function depend of the variables
$\frac{T_{m+1}}{\xi},\ldots,\frac{T_k}{\xi}$, then it means that
$\sigma_r$ is a total differential of some function depending on $x_1,x_2,\ldots,x_k$.

\vskip2mm
\subsection{7.1}
\vskip2mm

In this subsection we indicate a set of algebras in which the vectors
(\ref{e_1_e_2_e_3-k}) chosen arbitrarily and the equality (\ref{lambda-})
holds. We remind that an arbitrary commutative associative algebra, $\mathbb{A}_n^m$,
 with unit over the field of complex number $\mathbb{C}$ can be represented as
 $\mathbb{A}_n^m=S\oplus_s N$, where $S$ is $m$-dimensional semi-simple subalgebra
  and $N$ is $(n-m)$-dimensional nilpotent subalgebra (see Section~2).
\vskip2mm

\textbf{Theorem 5.}\label{teo-pro-napivprostu-alg-}
\textit{Let $\mathbb{A}_n^m=S\oplus_s N$. Then
the equality \em  (\ref{lambda-}) \em holds if at least one of the following conditions is satisfied:
\begin{enumerate}
    \item $\mathbb{A}_n^m\equiv S$;
    \item $N$ is a zero nilpotent subalgebra;
    \item $\dim_{\mathbb{C}}N\leq3$;
    \item $\dim_{\mathbb{C}}N=4$ and \end{enumerate}
 \begin{equation}\label{v-teor-pro-sigma-m+4-}
  \begin{array}{l}
\Upsilon_{m+1,m+2}^{m+1}\Upsilon_{m+2,m+3}^{m+2}=
\Upsilon_{m+1,m+2}^{m+1}\Upsilon_{m+2,m+4}^{m+2}=
\Upsilon_{m+1,m+3}^{m+1}\Upsilon_{m+3,m+4}^{m+2}=\vspace*{3mm} \\ 
=\Upsilon_{m+3,m+4}^{m+3}\Upsilon_{m+1,m+3}^{m+1}=
\Upsilon_{m+2,m+3}^{m+1}\Upsilon_{m+3,m+4}^{m+1}=
\Upsilon_{m+2,m+3}^{m+1}\Upsilon_{m+3,m+4}^{m+2}\vspace*{3mm} \\
=\Upsilon_{m+2,m+3}^{m+1}\Upsilon_{m+3,m+4}^{m+3}=
\Upsilon_{m+1,m+2}^{m+1}\Upsilon_{m+2,m+3}^{m+1}
\Upsilon_{m+3,m+4}^{m+2}= \vspace*{3mm} \\
=\Upsilon_{m+1,m+2}^{m+1}\Upsilon_{m+2,m+3}^{m+1}
\Upsilon_{m+3,m+4}^{m+3}
=\Upsilon_{m+2,m+3}^{m+2}\Upsilon_{m+3,m+4}^{m+1}=\vspace*{3mm} \\
=\Upsilon_{m+2,m+3}^{m+2}\Upsilon_{m+3,m+4}^{m+3}
=\Upsilon_{m+2,m+3}^{m+2}\Upsilon_{m+1,m+2}^{m+1}
\Upsilon_{m+3,m+4}^{m+1}=\vspace*{3mm} \\
=\Upsilon_{m+2,m+3}^{m+2}\Upsilon_{m+1,m+2}^{m+1}
\Upsilon_{m+3,m+4}^{m+2}
=\Upsilon_{m+2,m+3}^{m+2}\Upsilon_{m+1,m+2}^{m+1}
\Upsilon_{m+3,m+4}^{m+3}=0.\\
\end{array}
\end{equation}
  }\vskip 1mm

The proof is analogous to the proofs of Theorems 5~--- 8 of [\ref{Shpakivskyi-AACA-2015}].

Further we consider some examples of algebras, which satisfy the relations (\ref{v-teor-pro-sigma-m+4-}).\vskip2mm

\textbf{Examples.}
\begin{itemize}
  \item Consider the algebra with the basis $\{I_1:=1,I_2,I_3,I_4,I_5\}$ and
  multiplication rules:
  $$I_2^2=I_3\,,\,\,I_2\,I_4=I_5$$ and other products are zeros (for nilpotent
  subalgebra see [\ref{Martin}], Table 21, algebra $\mathcal{J}_{69}$ and
  [\ref{Burde_Fialowski}], page 590, algebra $A_{1,4}$).
  \item Consider the algebra with the basis $\{I_1:=1,I_2,I_3,I_4,I_5\}$ and
   multiplication rules:
  $$I_2^2=I_3$$ and other products are zeros (for nilpotent
  subalgebra see [\ref{Burde_Fialowski}], page 590, algebra $A_{1,2}\oplus A_{0,1}^2$).
  \item The algebra with the basis $\{I_1:=1, I_2, I_3, I_4, I_5\}$ and
   multiplication rules:
  $$I_2^2=I_3\,,\,\,I_4^2=I_5$$ and other products are zeros (for nilpotent
  subalgebra see [\ref{Burde_Fialowski}], page 590, algebra $A_{1,2}\oplus A_{1,2}$).
  \item The algebra with the basis $\{I_1:=1, I_2, I_3, I_4, I_5\}$ and
  multiplication rules:
  $$I_2^2=I_3\,,\,\,I_2\,I_3=I_4$$ and other products are zeros (for nilpotent
  subalgebra see [\ref{Martin}], Table 21, algebra $\mathcal{J}_{71}$).
\end{itemize}

In the paper [\ref{Shpakivskyi-AACA-2015}] is considered an example of algebra, which does not satisfy the relations (\ref{v-teor-pro-sigma-m+4-}). Moreover, in [\ref{Shpakivskyi-AACA-2015}] is selected the vectors $e_1,e_2,e_3$ of the form (\ref{e_1_e_2_e_3-k}) such that the equality
 (\ref{lambda-}) is not true.

 \vskip2mm
\subsection{7.2}
\vskip2mm

In this subsection we indicate sufficient conditions on a choose of
the vectors (\ref{e_1_e_2_e_3-k}) for which the equality (\ref{lambda-}) is true.
Let the algebra $\mathbb{A}_n^m$ be represented as $\mathbb{A}_n^m=S\oplus_s N$. Let us
note that the condition  $\zeta\in E_k\subset S$ means that in
the decomposition (\ref{e_1_e_2_e_3-k}) \,
 $a_{jr}=0$ for all $j=2,3,\ldots,k$ and $r=m+1,\ldots,n$.
\vskip2mm

\textbf{Theorem 6.}\label{teo-pro-napivprostu-alg-}
\textit{If $\mathbb{A}_n^m=S\oplus_s N$ and $\zeta\in E_k\subset S$,
 then the equality \em (\ref{lambda-}) \em holds.}\vskip 1mm

\textbf{\textit{Proof.}} Since $\zeta\in S$, then $T_s=0$ for $s=m+1,\ldots,n$ (see denotation
(\ref{B-})). We note that from the relation (\ref{obr-elem-0}) follows the equalities
\begin{equation}\label{sigma-k-}
\begin{array}{c}
\displaystyle
\sigma_{m+1}=\frac{dT_{m+1}}{\xi_{u_{m+1}}}+\widetilde{A}_{m+1}\,d\xi_{u_{m+1}}\,,\vspace*{4mm}\\ \displaystyle
\sigma_{r}=\frac{dT_r}{\xi_{u_r}}+\widetilde{A}_r\,d\xi_{u_r}+
\sum\limits_{q,s=m+1}^{r-1}\widetilde{A}_q\,dT_s\Upsilon_{q,r}^s\,,\quad r=m+2,\ldots,n.\\
\end{array}
\end{equation}

Now, from (\ref{lem_3_2---}), (\ref{A__p-}) follows that $\widetilde{A}_s=0$,
 and then from the (\ref{sigma-k-})  follows that $\sigma_r=0$ for $r=m+1,\ldots,n$. The equality (\ref{lambda-}) is a consequence of the equality $\sigma_r=0$ and the relation (\ref{obr-elem--1}).
\noindent The theorem is proved.

Let us note that by essentially the Theorem \ref{teo-pro-napivprostu-alg-}
generalizes the Theorem 3 of the paper [\ref{Shpakiv-Kuzm-Analele}] and generalizes the Theorem 9 of [\ref{Shpakivskyi-AACA-2015}].

Now we consider a case where $\zeta\notin S$. If $\mathbb{A}_n^m=S\oplus_s N$ and $\dim_{\mathbb{C}}N\leq3$,  then by Theorem  \ref{teo-pro-napivprostu-alg-}
 the equality (\ref{lambda-}) holds for any $\zeta\in E_k$.
 \vskip2mm

\textbf{Theorem 7.}\label{teo-pro-nil'potentnu-pidalg-rivne-4-dzeta}
\textit{Let $\mathbb{A}_n^m=S\oplus_s N$ and $\dim_{\mathbb{C}}N=4$.
Then the equality \em (\ref{lambda-}) \em holds if the following two conditions
satisfied:
  \begin{enumerate}
    \item $a_{j,m+1}=0$ for all $j=2,3,\ldots,k$;
    \item at least one of the relations $a_{j,m+2}=0$ or $a_{j,m+3}=0$ is true for all $j=2,3,\ldots,k$.
  \end{enumerate}}
 \vskip 1mm

The proof is similar to the proof of Theorem 10 of [\ref{Shpakivskyi-AACA-2015}].

\vskip 3.5mm

\section*{References}

\begin{enumerate}

\item \label{Shpakivskyi-AACA-2015}
	{\it Shpakivskyi V. S.} Curvilinear integral theorems for monogenic functions in commutative associative algebras // submitted to Adv. Appl. Clifford Alg., http://arxiv.org/pdf/1503.03464v1.pdf

\item \label{Lorch}
{\it Lorch E. R.}   The theory of analytic function in normed abelin
vector rings // Trans. Amer. Math. Soc., {\bf 54} (1943),
414~-- 425.

\item \label{Blum}
 {\it Blum E. K.} A theory of analytic functions in banach algebras //
 Trans. Amer. Math. Soc., {\bf 78} (1955), 343~-- 370.

\item \label{Pl-Shp3}
{\it Shpakivskyi V. S., Plaksa S. A.}  Integral theorems and a
Cauchy formula in a commutative three-dimensional harmonic
algebra // Bulletin Soc. Sci. Lettr. L\'od\'z, {\bf 60} (2010),
47~-- 54.

\item \label{Pl-Shp-Algeria}
{\it Plaksa S. A., Shpakivskyi V. S.}  Monogenic functions in a
finite-dimensional algebra with unit and radical of maximal
dimensionality // J. Algerian Math. Soc., {\bf 1} (2014), 1~-- 13.

\item \label{Pl-Pukh-Analele}
{\it Plaksa S. A., Pukhtaievych R. P.}   Constructive
description of monogenic functions in $n$-dimensional semi-simple
algebra //  An. \c{S}t. Univ. Ovidius Constan\c{t}a, \textbf{22} (2014),
 no.~1, 221~-- 235.

\item \label{Ketchum-28}
\textit{Ketchum P. W.}  Analytic functions of hypercomplex variables //
Trans. Amer. Math. Soc., {\bf 30} (1928), no. 4, 641~-- 667.

\item \label{Ketchum-29}
\textit{Ketchum P. W.} A complete solution of Laplace's equation
by an infinite hypervariable // Amer. J. Math., \textbf{51} (1929),
179~-- 188.

\item \label{Rosculet-54}
\textit{Ro\c{s}cule\c{t} M. N.}  O teorie a func\c{t}iilor de o variabil\u{a}
hipercomplex\u{a} \^{i}n spa\c{t}iul cu trei dimensiuni // Studii \c{s}i
Cercet\v{a}ri Matematice, \textbf{5}, nr. 3--4 (1954), 361~-- 401.

\item\label{Rosculet-55}
\textit{Ro\c{s}cule\c{t} M. N.} Algebre liniare asociative \c{s}i  comutative
\c{s}i finc\c{t}ii monogene ata\c{s}ate lor // Studii \c{s}i Cercet\v{a}ri
Matematice, \textbf{6}, nr. 1--2 (1955), 135~-- 173.

\item\label{Cartan}
\textit{Cartan E.} Les groupes bilin\'{e}ares et les syst\`{e}mes de nombres complexes //
Annales de la facult\'{e} des sciences de Toulouse, \textbf{12} (1898), no.~1, 1~-- 64.

\item \label{Burde_de_Graaf}
\textit{Burde D., de Graaf W.} Classification of Novicov algebras // Applicable Algebra
in Engineering, Communication and Computing, \textbf{24}(2013), no.~1, 1~-- 15.

\item \label{Burde_Fialowski}
\textit{Burde D., Fialowski A.} Jacobi–Jordan algebras //
Linear Algebra Appl., \textbf{459} (2014), 586~-- 594.

\item \label{Martin}
\textit{Martin M. E.} Four-dimensional Jordan algebras //
Int. J. Math. Game Theory Algebra \textbf{20} (4) (2013)
41~-- 59.

\item \label{Shpakivskyi-2014}
\textit{Shpakivskyi V. S.}  Constructive description of monogenic functions in
a finite-dimensional commutative associative algebra // submitted to J. Math. Anal. Appl.,\, http://arxiv.org/pdf/1411.4643v1.pdf

\item \label{Hil_Filips}
 {\it Hille E.,  Phillips R. S.} Functional analysis and semi-groups
 [Russian translation],
Inostr. Lit., Moscow (1962).

\item \label{Shpakivskyi-Zb-2015-1}	
 {\it Shpakivskyi V. S.} Monogenic functions in finite-dimensional commutative associative algebras // accepted to Zb. Pr. Inst. Mat. NAN Ukr.

 \item \label{Privalov}
{\it Privalov I. I.}  Introduction to the Theory of Functions of a Complex Variable,
GITTL, Moscow, 1977. (Russian)

\item \label{Shabat}
{\it Shabat B. V.} Introduction to complex analysis, Part 2,
Nauka, Moskow (1976). (Russian)

\item\label{Pl-Pukh}
{\it Plaksa S. A., Pukhtaevich R. P.}  Constructive
description of monogenic functions in a three-dimensional harmonic
algebra with one-dimensional radical // Ukr. Math. J., {\bf 65}
(2013), no.~5, 740~-- 751.

\item \label{Shpakiv-Kuzm-Analele}
{\it Shpakivskyi V. S., Kuzmenko T. S.}  Integral theorems for the quaternionic
 $G$-monogenic mappings // accepted to An. \c{S}t. Univ. Ovidius Constan\c{t}a.

\end{enumerate}
\label{end}
\end{document}